\documentclass[twocolumn,preprintnumbers,amsmath,amssymb]{revtex4}

\usepackage{graphicx}% Include figure files
\usepackage{dcolumn}% Align table columns on decimal point
\usepackage{bm}% bold math

\begin{document}

\title{
One class of conservative difference schemes for solving molecular dynamics equations of motion}

%%\shorttitle{Surface Tension of QG Bag and the Colour Tube String Tension} 

\author{E. G. Nikonov}

\affiliation{Laboratory of Information Technologies, Joint Institute for Nuclear Research, Dubna, Russia}

\vspace{1.cm}

\begin{abstract}
Simulation of many-particle system evolution by molecular dynamics takes to decrease integration step to provide numerical scheme stability on the sufficiently large time interval. It leads to a significant increase of the volume of calculations. An approach for constructing symmetric simplectic numerical schemes with given approximation accuracy in relation to integration step, for solving molecular dynamics Hamiltonian equations, is proposed in this paper. Numerical experiments show that obtained under this approach symmetric simplectic third order scheme is more stable for integration step, time-reversible and conserves Hamiltonian of the system with more accuracy at a large integration interval then second order velocity Verlet numerical schemes. 
\\

\end{abstract}

\maketitle

%%%%%%%%%%%%%%%%%%%%%  Introduction

%%%%%%     04.05.2016        Eduard G. Nikonov             SSDSHEq.tex

%%%%%%%%%%%%%%%%%%

\section{Introduction}

Numerical schemes, which is using for solving systems of equations of many-particle dynamics, can have restrictions on a step and an interval of integration because if they increase, the numerical schemes become unstable and do not conserve integrals of motion. As a result, when we simulate many-particle system behaviour on the sufficiently large time interval we should decrease an integration step, which leads to considerable increasing of computation quantity.  

\section{HAMILTONIAN MOLECULAR DYNAMICS}
A motion of $N$-particle system in a field with potential $V(\mathbf{q})$ can be described by the system of Hamiltonian equations \citep{Hairer:06}
\begin{equation}\label{Heq}
\dot{\mathbf{p}} = -\frac{\partial H(\mathbf{p},\mathbf{q})}{\partial \mathbf{q}},
\,\,
\dot{\mathbf{q}} = \frac{\partial H(\mathbf{p},\mathbf{q})}{\partial \mathbf{p}},
\end{equation}
with initial conditions $\mathbf{p}(0)=\mathbf{p}^0,\,\mathbf{q}(0)=\mathbf{q}^0$.

Here $\mathbf{q}=(q_1,…,q_d)^T$- particle coordinates, $\mathbf{p}=(p_1,… ,p_d)^T$- particle momentums, $d=3N$ - a dimension of coordinate space and
$$ 
H(\mathbf{p},\mathbf{q})=\frac{1}{2}\mathbf{p}^T M^{-1}(\mathbf{q})\ \mathbf{p}+V(\mathbf{q}) 
$$ 
is a separable Hamiltonian of the system with symmetric and positive definite mass matrix $M(\mathbf{q})$ and a field potential $V(\mathbf{q})$.

\section{DESIGN OF CONSERVATIVE DIFFERENCE SCHEMES}
An approach to design of conservative difference scheme for solving Hamiltonian equations (\ref{Heq}) is based on following stages. First is a choice of appropriate type of generating function $S$, which fixes a definite family of simplectic difference scheme \cite{Kang:1989}. Simplectic difference scheme corresponds to the canonical transformation of canonical variables \cite{Gantmakher}. Solution of Hamiltonian system (\ref{Heq}) for a definite time moment can be represented as canonical transformation and that is why conserves the value of Hamiltonian \cite{Arnold}. Second is using of “forward” and “backward” Taylor expansion of coordinate and momentum on a time step for getting of a corresponding generating function and remaining values of coordinate and momentum. At last, obtained at the previous stage explicit and implicit schemes are used for constructing symmetric simplectic scheme for solving Hamiltonian equations (\ref{Heq}).
$$
S  =  \frac{h}{2}\left(\frac{\mathbf{q}^{k+1}-\mathbf{q}^k}{h}\right)^2 - \frac{h}{2}\left[\ V(\mathbf{q}^k) + V(\mathbf{q}^{k+1})\ \right] -
$$
$$ - \frac{h}{12}\left[\ \nabla V(\mathbf{q}^{k+1})-\nabla V(\mathbf{q}^k)\ \right] \cdot\left(\mathbf{q}^{k+1}-\mathbf{q}^k\right),
$$

$$
\mathbf{p}^k   =  \frac{\mathbf{q}^{k+1}-\mathbf{q}^k}{h}-\frac{h}{12}\left[\ 5\nabla V(\mathbf{q}^k)+\nabla V(\mathbf{q}^{k+1})\ \right]- 
$$
$$ 
-\frac{h}{12}\mathcal{H}(V)(\mathbf{q}^k)\cdot(\mathbf{q}^{k+1}-\mathbf{q}^k),
$$
$$
\mathbf{p}^{k+1}   =   \frac{\mathbf{q}^{k+1}-\mathbf{q}^k}{h}+\frac{h}{12}\left[\ \nabla V(\mathbf{q}^k)+5\nabla V(\mathbf{q}^{k+1})\ \right] -
$$
$$
- \frac{h}{12}\mathcal{H}(V)(\mathbf{q}^{k+1})\cdot(\mathbf{q}^{k+1}-\mathbf{q}^k),
$$
where $\mathcal{H}(V)(q)$ is Hessian matrix for $V(q)$.
 
\section{NUMERICAL EXPERIMENTS}
Obtained by described above approach difference scheme of third order [5] was tested by solving the system of equations (\ref{Heq}) with Hamiltonian
$$
H(\mathbf{p},\mathbf{q})=\frac{1}{2}(p_1^2+p_2^2)-\frac{1}{\sqrt{q_1^2+q_2^2}},
$$ 
so called Kepler two body problem,
Results of numerical calculations are compared for the symmetric simplectic scheme of third order and well-known second order velocity Verlet scheme with a time step $h=0.2$ and a time interval $[0,T], \,T=5000$. Approximate solution by Verlet scheme with a time step h=0.02 is used as an exact solution [5].
\begin{figure}[h]
\center{\includegraphics[width=1\linewidth]{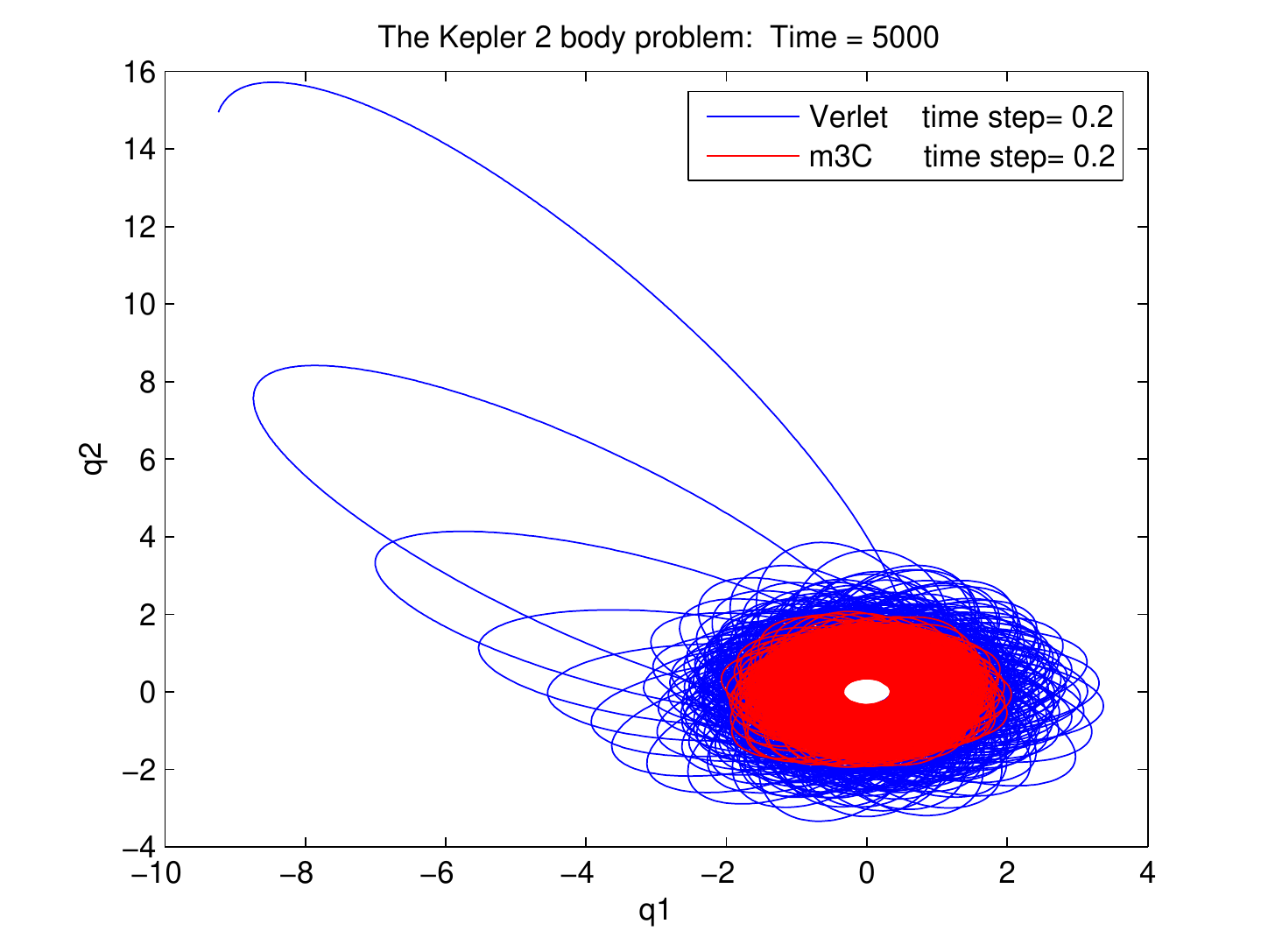}} 
\caption{Phase trajectories.}
\label{Traject}
\end{figure}
\begin{figure}[h] 
\center{\includegraphics[width=1\linewidth]{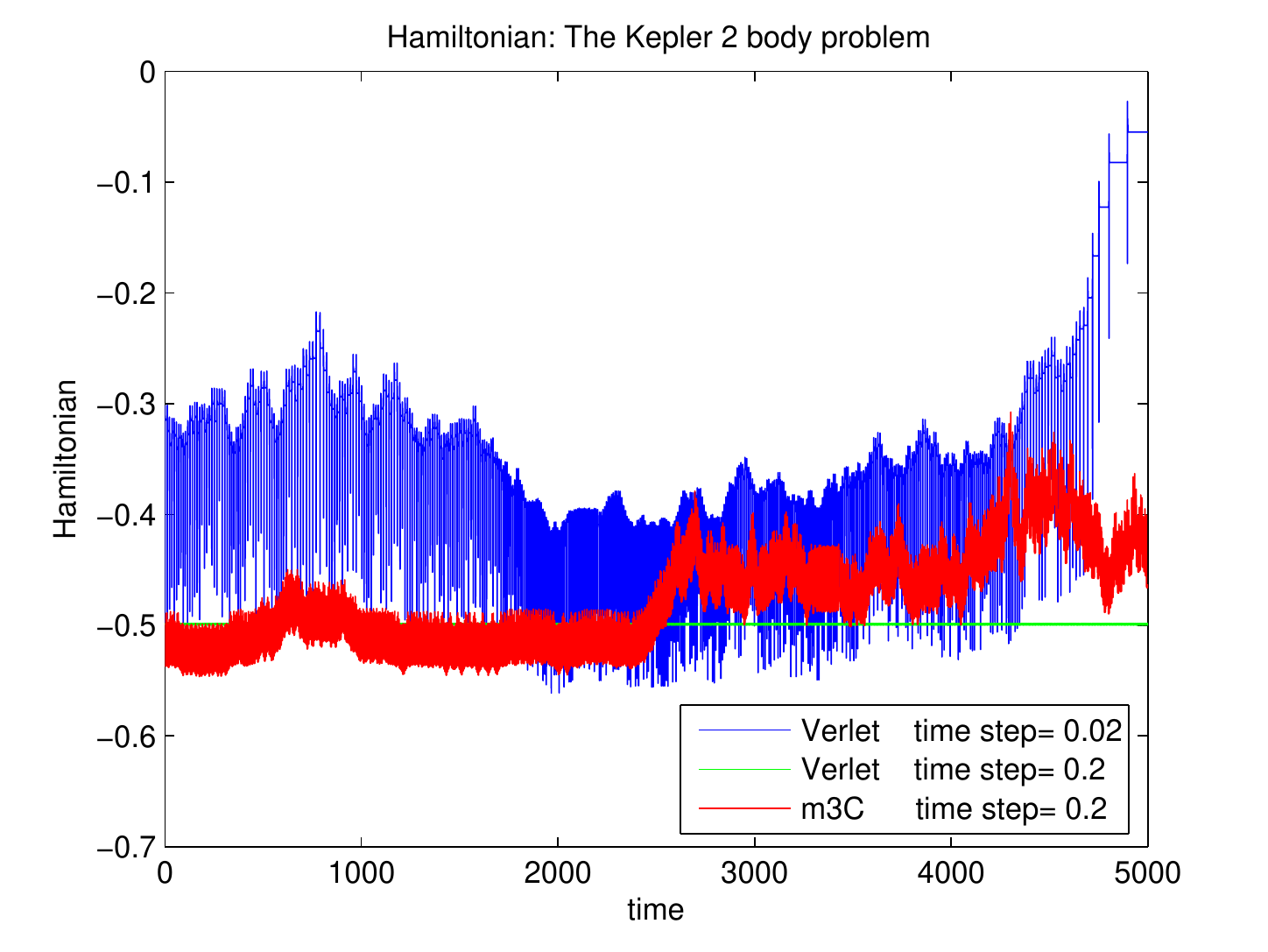}} 
\caption{Conservation of Hamiltonian.}
\label{Humiltonian}
\end{figure} 

\section{Results}

A new approach for constructing symmetric simplectic numerical schemes for solving Hamiltonian systems of equations is proposed. The numerical schemes constructed by this approach produce more stable and accurate solution of Hamiltonian system (1) and better conserve the energy of a system on the large interval of numerical integration for a relatively large integration step in comparision with the Verlet method, which is usually using for solving equations of motion in molecular dynamics. 
%\vskip3mm

{\bf Acknowledgments.} E.G.N. acknowledge Dr. B.~Batgerel for computations and numerical results. The work was supported in part by the Russian
Foundation for Basic Research, Grant No. 15-01-06055.

%\vskip3mm

\end{document}